\documentclass{amsart}

\usepackage{amssymb}

\newtheorem{thm}{Theorem}[section]

\newtheorem{lem}[thm]{Lemma}
\newtheorem{cor}[thm]{Corollary}
\theoremstyle{definition}

\theoremstyle{remark}

\numberwithin{equation}{section}

\begin{document}

\title{Large time behavior of solutions to nonlinear beam equations}

\author{David Raske}
\email{nonlinear.problem.solver@gmail.com}

\begin{abstract}

In this article we will investigate the large time behavior of solutions of a special class of initial/boundary value problems that involve nonlinear damped beam equations. We will show that the solution energies of global pseudo classical solutions to these initial/boundary value problems decay exponentially.

\end{abstract}                        

\maketitle

\section{Introduction}
Let $a$, $b$ be two real numbers. Let $F$ and $G$ be two functions from $\mathbb{R}$ into $\mathbb{R}$. Let $f$ be           a function from $(a,b) \times(0,\infty)$ into $\mathbb{R}$. Initial/boundary value problems of the form
\begin{equation}\label{pseudo}
\begin{split}
& u_{tt} + u_{xxxx} + F(u_t) + G(u) = f(x,t), \text{ for all } (x,t) \in (a,b) \times (0,\infty) \\
& u(a,t)=0=u(b,t),  \text{ for all } t \text{ in } (0,\infty) \\
&  u_{xx}(a,t)=0=u_{xx}(b,t),  \text{ for all } t \text{ in } (0,\infty) \\
& u(x,0)=u_0(x), u_t(x,0) = u_1(x) \text{ for all } x \text{ in } (a,b). \\
\end{split}
\end{equation}
arise naturally in the study of vibrations in suspension bridges. The issue of whether or not the above initial/boundary value problem is well-posed is addressed in \cite{R}.

Before we begin our investigation of the large time behavior of solutions to the above initial/boundary problem, we need to define some spaces. Let $\Omega$ be the open, bounded domain in $\mathbb{R}$, $(a,b)$. Let $k$ be a positive integer and let $H^k(\Omega)$ consist of all locally summable functions $u: \Omega \rightarrow \mathbb{R}$ such that for each multindex $\alpha$ with $|a| \leq k$, $D^\alpha u$ exists in the weak sense and belongs to $L^2(\Omega)$. Let $H^k_0(\Omega)$ be the closure of $C^\infty_c(\Omega)$ in $H^k(\Omega)$. Now let, $H^2_*(\Omega)$ be the intersection of $H^1_0(\Omega)$ with $H^2(\Omega)$. We will equip it with the inner-product  
\begin{equation}\label{first}
(u,v)_{H^2_*} = \int_\Omega u_{xx} v_{xx} dx.
\end{equation}
Let $H^4_*(\Omega)$ be the elements $u$ of $H^4(\Omega)$ such that $u \in H^2_*(\Omega)$ and $u_{xx} \in H^2_*(\Omega)$; we will equip it with the inner-product 
\begin{equation}\label{second}
(u,v)_{H^4_*} = \int_\Omega u_{xxxx} v_{xxxx} dx.
\end{equation}
 In \cite{R} it is shown that the inner-product \eqref{first} makes $H^2_*(\Omega)$ a Hilbert space and the inner-product \eqref{second} makes $H^4_*(\Omega)$ a Hilbert space.


Now, suppose there exists a function $\mathbf{u}: [0,\infty) \rightarrow L^2((a,b))$ such that $\mathbf{u} \in C( [0,\infty) ,H^4_*((a,b)) )$, $\mathbf{u} \in C^1( [0,\infty) ,H^2_*((a,b)) )$, $\mathbf{u} \in C^2( [0,\infty), L^2((a,b)) )$, $\mathbf{u}(0)=u_0$ with $u_0 \in H^4_*((a,b))$, $\mathbf{u}'(0)=u_1$ with $u_1 \in H^2_*((a,b))$, and such that for all $\phi \in C^\infty_c(0,\infty)$ and $v \in L^2((a,b))$ we have
\begin{equation}\label{weak}
\begin{split}
&\int_0^\infty (\mathbf{u}''(t),v)_{L^2}\phi(t) \, dt \\
& = -\int_0^\infty [(L(\mathbf{u}(t)),v)_{L^2} + (F(\mathbf{u'}(t)),v)_{L^2} + (G(\mathbf{u}(t)),v)_{L^2} - ( \mathbf{f} (t),v)_{L^2}]\phi(t) \, dt.
\end{split}
\end{equation}
Here $L$ is the $L^2((a,b))$-valued differential operator $\Delta^2$, where $\Delta$ is the Laplacian, whose domain is $H^{4}_*((a,b))$.  $\mathbf{f}(t) \in C([0,\infty); L^2((a,b))).$ We will then call $\mathbf{u}(t)$ a global pseudo classical solution of the initial/boundary value problem \eqref{pseudo}. Note that if $\mathbf{u}(t)$ is a pseudo classical solution of \eqref{pseudo} we have it that
\begin{equation}\label{weak2}
 (\mathbf{u}''(t),v)_{L^2} +  (F(\mathbf{u}'(t)),v)_{L^2} + (L(\mathbf{u}(t)),v)_{L^2}  = (\mathbf{f}(t),v)_{L^2}
\end{equation}
for each $v \in L^2((a,b))$ and all $t \in (0,\infty)$.


If a global pseudo classical solution to the initial/boundary value problem \eqref{pseudo} exists, a natural question to ask is whether or not it converges to a solution of the boundary value problem
\begin{equation} \label{eq 1.1}
\begin{split}
& u_{xxxx} + G(u) = f(x) \text{ on } (a,b),  \\
& u(a)=0=u(b), \\
& u_{xx}(a)=0=u_{xx}(b). \\
\end{split}
\end{equation}
as $t \rightarrow \infty$. Since the source of motivation to study problems of the form \eqref{pseudo} oftentimes comes from continuum mechanics, another important question to ask is whether or not the energy of a solution to \eqref{pseudo} is decreasing if $\mathbf{f}(t)$ is independent of time. These two questions will be investigated in this paper. Once we are done, we will have 

\begin{thm} Let $c$ and $d$ be two real numbers, with $c < d$. Let $m$ be a positive real number such that $m \geq 2$. Let $a_1$,and $a_2$ be two  positive real numbers with $a_1 \leq a_2$. Let $U$ be the open interval $(c,d)$. Let $\mathbf{f}(t)$ be a member of $C([0,\infty); L^2((a,b)))$, such that $\mathbf{f}(t)  \equiv 0$ for all $t \in [0,\infty)$. Let $F: \mathbb{R} \rightarrow \mathbb{R}$ be a continuous function such that $F(z)z\geq 0$ for all $z \in \mathbb{R}$ and such that
$$
|a_1 (z + |z|^{m-2}z)| \leq |F(z)| \leq |a_2 (z+|z|^{m-2}z)|.
$$
for all $z \in \mathbb{R}$. Suppose $G: \mathbb{R} \rightarrow \mathbb{R}$ has the property $G \equiv 0$. Furthermore, let $u_0$ be an element of $H_*^4(U)$ and let $u_1$ be an element of $H^2_*(U)$. Then, if there exists a unique global pseudo classical solution, $\mathbf{u}(t)$, of the initial/boundary value problem \eqref{pseudo}, we have that $||\mathbf{u}(t)||_{H^2_*((a,b))}$ exhibits exponential decay and $||\mathbf{u}'(t)||_{L^2((a,b))}$ exhibits exponential decay.
\end{thm}

An immediate consequence of the above is the following

\begin{cor} Let $c$ and $d$ be two real numbers, with $c < d$. Let $m$ be a positive real number such that $m \geq 2$. Let $a_1$,and $a_2$ be two  positive real numbers with $a_1 \leq a_2$. Let $U$ be the open interval $(c,d)$. Let $\mathbf{f}(t)$ be a member of $C([0,\infty);L^2((a,b)))$ such that $\mathbf{f}(t) \equiv  \mathbf{f}(0)$ for all $t \in [0,\infty)$. Let $F: \mathbb{R} \rightarrow \mathbb{R}$ be a continuous function such that $F(z)z \geq 0$ for all $z \in \mathbb{R}$ and such that
$$
|a_1 (z + |z|^{m-2}z)| \leq |F(z)| \leq |a_2 (z+|z|^{m-2}z)|.
$$
for all $z \in \mathbb{R}$. Suppose $G: \mathbb{R} \rightarrow \mathbb{R}$ has the property $G \equiv 0$. Furthermore, let $u_0$ be an element of $H_*^4(U)$ and let $u_1$ be an element of $H^2_*(U)$. Then, if there exists a unique global pseudo classical solution, $\mathbf{u}(t)$, of the initial/boundary value problem \eqref{pseudo}, we have that $||\mathbf{u}(t)-\hat{u}||_{H^2_*((a,b))}$ exhibits exponential decay, where $\hat{u}$ is the solution of \eqref{eq 1.1}, and $||\mathbf{u}'(t)||_{L^2((a,b))}$ exhibits exponential decay.
\end{cor}

In section two of this paper we make some observations about the energy of solutions to \eqref{pseudo}. In particular, we will see that it is nonincreasing. A consequence of this is that the $C(\overline{U})$ norm of the solution is bounded independent of $t \in [0,\infty)$. We will also prove an inequality that will become very useful when combined with the boundedness of solutions.

In section three of this paper we will prove Theorem 1.1. The idea behind the proof was borrowed from \cite{M}. One perturbs the energy by a small quantity and then shows that this modified energy decays exponentially. The only results needed are the standard inequalities one uses while studying partial differential equations along with the inequality that is proved in section two. The main difference between the proof presented below and the one found in \cite{M} is that we do not have to consider different cases of the size of the $L^m$ norm of the solutions as $t$ increases. This is because we do not consider equations with source terms.

\section{Preliminaries}

\subsection{The Solution Energy}
Let $\mathbf{u}(t)$ be a global pseudo classical solution to \eqref{pseudo}, where $F$ satisfies the hypotheses contained in the statement of Theorem 1.1, where $G \equiv 0$, and where $\mathbf{f}(t) \equiv 0$ for all $t \in [0,\infty)$. Set $[\mathbf{u}(t)](x) = u(x,t) (x \in (a,b), 0 \leq t)$.
We will call 
\begin{equation}
E(t) = \frac{1}{2}\int_U (u_t(x,t))^2 \, dx + \frac{1}{2} \int_U (u_{xx}(x,t))^2 \, dx,
\end{equation}
the solution energy. In important property of the solution energy is that
\begin{equation}
E'(t) = -\int_U F(u_t)u_t(x,t) \, dx \leq 0.
\end{equation}
An important consequence of this property is that there exists a positive real number $C$ such that $||\mathbf{u}(t)||_{C(\overline{U})} \leq C$ for all $t \in [0,\infty)$.

\subsection{Inequalities}

An inequality that will be used multiple times in the proof of Theorem 1.1 is as follows. There exists a positive real number $B$ such that
\begin{equation}\label{Poincare}
||u||_2 \leq B ||u||_{H^2_*},
\end{equation}
for all $u \in H^2_*(U)$. \footnote{Let $v \in L^2((a,b))$. Here, and for the remainder of the paper, $||v||_2$, will denote the $L^2$ norm of $v$.} This inequality is proven in Section Two of \cite{R}.

Another inequality that will prove to be useful below is the following

\begin{lem}
Let $n$ be a positive integer, and let $m$ be a positive real number such that $m \geq 2$. Let $U$ be a bounded, open subset of $\mathbb{R}^n$ with $C^1$ boundary. Let $z$ be an element of $C([0,\infty); C(\overline{U}))$ such that there exists a positive real number $M_1$ such that $||z(t)||_{C(\overline{U})} \leq M_1$ for all $t \in [0,\infty)$. Then we have the existence of a positive real number $C$  such that
$$
\int_U |z(t)|^m \, dx \leq C \int_U  (z(t))^2 \, dx,
$$
for all $t \in [0,\infty)$. 
\end{lem}

\begin{proof}
Let $v$ be an element of $C([0,\infty); C(\overline{U}))$ such that there exists a positive real number $M_2$ such that $||v(t)||_{C(\overline{U})} \leq M_2$ for all $t \in [0,\infty)$. Let $M:=\max \{M_1,M_2\}$. Let $h: \mathbb{R} \rightarrow \mathbb{R}$ be a function such that $h(x) = |x|^{m/2}$. Since $h$ is locally Lipschitz we have the existence of a positive real number $L_h$ such that
$$
|h(s_1)-h(s_2)| \leq L_h|s_1 - s_2| \text{ for any } s_1, s_2 \in [-M,M].
$$
This, in turn, allows us to write
$$
||h(z(t)) - h(v(t))||_{2}^2 \leq ( L_h)^2 ||z(t) - v(t)||_{2}^2
$$
for all $t \in [0,\infty)$. Setting $v(t) \equiv 0$ for all $t \in [0,\infty)$, we obtain the lemma.
\end{proof}

\section{Proof of Theorem 1.1}
First, set $[\mathbf{u}(t)](x) = u(x,t)$. Proceeding as in \cite{M}, we start with the following observation. Notice that 
\begin{equation}\label{energy derivative}
\begin{split}
E'(t)&  = - \int_U F(u_t) u_t(x,t)  \, dx \\     
       & \leq - a_1( \int_U (u_t(x,t))^2 \, dx + \int_U |u_t(x,t)|^m \, dx).
\end{split}
\end{equation}
Define 
\begin{equation}\label{def}
H(t) = E(t) + \epsilon \int_U u u_t(x,t) \, dx
\end{equation}
for an $\epsilon$ to be specified later. Notice that we can assume without loss of generality that $B \geq 1$, where $B$ is the constant appearing in \eqref{Poincare}, so let us do so. Then the Schwarz inequality gives us
\begin{equation}\label{Schwarz}
\begin{split}
|\epsilon \int_U u u_t(x,t) \, dx| & \leq \frac{\epsilon}{2} (\int_U u(x,t)^2 \, dx) + \frac{\epsilon}{2} \int_U (u_t(x,t))^2 \, dx \\ & \leq \frac{\epsilon B^2}{2} \int_U (u_{xx}(x,t))^2 \, dx + \frac{\epsilon B^2}{2} \int_U (u_t(x,t))^2 \,dx \\ & \leq \epsilon B^2 E(t).
\end{split}
\end{equation}
It follows that
\begin{equation}\label{H-E}
|H(t) - E(t)| \leq \epsilon B^2 E(t).
\end{equation}
Next we differentiate \eqref{def} and use \eqref{weak2} and \eqref{energy derivative} to see that
\begin{equation}\label{ineq}
\begin{split}
H'(t) \leq & -a_1 \int_U (u_t(x,t))^2 \, dx - a_1 \int_U (u_t(x,t))^m \,dx + \epsilon \int_U (u_t(x,t))^2 \,dx \\ & - \epsilon \int_U (u_{xx}(x,t))^2 \,dx - \epsilon \int_U F(u_t)u(x,t) \,dx \\  & \leq - a_1 \int_U (u_t(x,t))^2 \,dx - a_1 \int_U |u_t(x,t)|^m \,dx + \frac{3}{2} \epsilon \int_U (u_t(x,t))^2 \,dx \\ &  - \frac{1}{2} \epsilon \int_U (u_{xx}(x,t))^2 \,dx -  \epsilon \int_U F(u_t)u (x,t) \,dx - \epsilon E(t) \\ & \leq  - a_1 \int_U (u_t(x,t))^2 \,dx - a_1 \int_U |u_t(x,t)|^m \,dx + \frac{3}{2} \epsilon \int_U (u_t(x,t))^2 \,dx \\ &  - \frac{1}{2} \epsilon \int_U (u_{xx}(x,t))^2 \,dx + \epsilon a_2 \int_U |u_t||u|(x,t) \,dx + \epsilon a_2 \int_U |u_t|^{m-1}|u|(x,t) \,dx \\ & - \epsilon E(t). 
\end{split}
\end{equation}
By using 
\begin{equation}\label{Cauchy}
a_2 \int_U |u_t||u|(x,t) \, dx \leq \frac{1}{4} \int_U (u_{xx}(x,t))^2 + a_2^2 B \int_U |u_t(x,t)|^2 \, dx,
\end{equation}
inequality \eqref{ineq} takes then the form
\begin{equation}\label{ineq 2}
\begin{split}
H'(t) & \leq -a_1 \int_U |u_t(x,t)|^2 \, dx - a_1 \int_U |u_t(x,t)|^m \, dx \\ & + (\frac{3}{2} + (a_2^2 B))\epsilon \int_U (u_t(x,t))^2 \, dx - \frac{1}{4} \epsilon \int_U (u_{xx}(x,t))^2 \,dx \\ & - a_2 \epsilon \int_U |u_t|^{m-1}|u|(x,t) \, dx -\epsilon E(t).
\end{split}
\end{equation}
We then exploit Young's inequality
\begin{equation}\label{Young's}
XY \leq \delta X^r + c(\delta)Y^s,
\end{equation}
where $X, Y, \delta, c(\delta) \geq 0$ and $\frac{1}{r} + \frac{1}{s} = 1$, with $r = m$ and $s = \frac{m}{m-1}$ to get
\begin{equation}\label{consequence}
\int_U |u_t|^{m-1}|u|(x,t) \,dx \leq \delta ||u||^m_m + c(\delta)||u_t||^m_m,
\end{equation}
for all $\delta > 0$. We can now combine the fact that there exists a positive real number $C$ such that $||\mathbf{u}(t)||_{C(\overline{U})} \leq C$ for all $t \in [0,\infty)$ with Lemma 2.1 to see that there exists a positive real number $\gamma$ such that
\begin{equation}\label{ineq 3}
\int_U |u_t|^{m-1}|u|(x,t) \,dx \leq \delta \gamma ||u||_2^2 + c(\delta) ||u_t||^m_m.
\end{equation}
Therefore \eqref{ineq 2} becomes
\begin{equation}\label{ineq 4}
\begin{split}
H'(t) & \leq -a_1 \int_U |u_t(x,t)|^2 \,dx - a_1 \int_U |u_t(x,t)|^m \,dx \\  & + (\frac{3}{2} + a^2_2 B) \epsilon \int_U (u_t(x,t))^2 \,dx - \frac{1}{4} \epsilon \int_U (u_{xx}(x,t))^2 \,dx \\  & + a_2 \epsilon (\gamma \delta ||u||^2_2 + c(\delta)||u_t||^m_m) - \epsilon E(t),
\end{split}
\end{equation}
for all $\delta > 0$. This, in turn, allows us to write
\begin{equation}\label{ineq 5}
\begin{split}
H'(t) \leq & -\epsilon E(t) - \frac{1}{4} \epsilon ||u_{xx}||^2_2 + a_2 \epsilon \gamma \delta ||u||^2_2 \\ &  - [a_1 - (\frac{3}{2} + a_2^2 B)\epsilon]||u_t||_2^2 - a_1[1- \frac{a_2}{a_1} \epsilon c(\delta)] ||u_t||_m^m,
\end{split}
\end{equation}
for all $\delta > 0$. Now, picking $\delta \leq \frac{B^2}{4 a_2 \gamma}$, we see that
\begin{equation}\label{ineq 6}
H'(t)  \leq -\epsilon E(t) - a_1[1 - \frac{a_2}{a_1} \epsilon c(\delta)]||u_t||_m^m - [a_1 - (\frac{3}{2} + a^2_2 B)\epsilon]||u_t||_2^2.
\end{equation} 
Now, if we pick $\epsilon \leq \min\{\frac{a_1}{a_2}\frac{1}{c(\delta)},\frac{a_1}{3/2 + a_2^2 B}\}$, we see that $H'(t) \leq - \epsilon E(t)$. Now recall that 
\begin{equation}\label{difference}
|H(t)-E(t)| \leq \epsilon B^2 E.
\end{equation}
It follows that $E(t) \geq \frac{1}{1+\epsilon B^2} H(t)$, and hence
\begin{equation}\label{H decay equation}
H'(t) \leq (-\epsilon \frac{1}{1+\epsilon B^2}) H(t).
\end{equation}
Calculus thus gives us that there exists a positive real number $r$ such that
\begin{equation}\label{H decay}
H(t) \leq H(0) e^{-rt}.
\end{equation}
Now, recalling \eqref{difference}, we see that
\begin{equation}\label{E bound}
E(t) \leq \frac{H(t)}{(-\epsilon B^2 + 1)},
\end{equation}
provided that $\epsilon$ is chosen small enough so that $-\epsilon B^2 + 1 > 0$. It follows that if we pick $\epsilon$ small enough, we have 
\begin{equation}\label{final}
E(t) \leq \frac{H(t)}{-\epsilon B^2 + 1} \leq \frac{H(0)}{-\epsilon B^2 + 1} e^{-rt},
\end{equation}
with $-\epsilon B^2 + 1 > 0$. The theorem follows.  

\bibliographystyle{amsplain}

\end{document}